\documentclass[reqno]{amsart}
\usepackage{amscd}
\usepackage[dvips]{graphics}
\usepackage{mathrsfs}
\usepackage{bm}
\usepackage[sort&compress,numbers]{natbib}  

\def\beq{\begin{equation}}
\def\eeq{\end{equation}}
\def\ba{\begin{array}}
\def\ea{\end{array}}

\numberwithin{equation}{section}
\newenvironment{abs}{\textbf{Abstract}\mbox{  }}{ }
\newenvironment{key words}{\textbf{Keywords}\mbox{  }}{ }
\newtheorem{theorem}{Theorem}[section]

\newtheorem{proposition}[theorem]{\textbf{Proposition}}

\theoremstyle{remark}
\newtheorem{remark}[theorem]{\textbf{Remark}}
\theoremstyle{plain}

\allowdisplaybreaks[4]

\begin{document}
\title[Sharp Sobolev inequalities on the complex sphere]{Sharp Sobolev inequalities on the complex sphere}
\author  {Yazhou Han and Shutao Zhang}
\address{Yazhou Han, Department of Mathematics, College of Science, China Jiliang University, Hangzhou, 310018, China}

\email{yazhou.han@gmail.com}

\address{Shutao Zhang, Department of Mathematics, College of Science, China Jiliang University, Hangzhou, 310018, China}

\email{taoer558@163.com}

\date{}
\maketitle

\noindent
\begin{abs}
This paper is devoted to establish a class of sharp Sobolev inequalities on the unit complex sphere
as follows:

1) {\bf Case $0<d<Q=2n+2$}: for any $f\in C^\infty$ and $2\leq q \leq \frac{2Q}{Q-d}$,
\begin{align*}
    \|f\|_q^2\leq & \frac{8(q-2)}{d(Q-d)} \frac{\Gamma^2((Q-d)/4+1)} {\Gamma^2((Q+d)/4)}\left( \int_{\mathbb{S}^{2n+1}} f\mathcal{A}_df d\xi\right. \nonumber\\
    &\left. -\frac{\Gamma^2((Q+d)/4)} {\Gamma^2((Q-d)/4)} \int_{\mathbb{S}^{2n+1}} |f|^2 d\xi\right) +\int_{\mathbb{S}^{2n+1}} |f|^2 d\xi;
\end{align*}

2) {\bf Case $d=Q$}: for any $f\in C^\infty \cap\mathbb{R}\mathcal{P}$ and $2\leq q< +\infty$,
\begin{equation*}
    \|f\|_q^2\leq \frac{q-2}{(n+1)!} \int_{\mathbb{S}^{2n+1}} f \mathcal{A}'_Q f d\xi +\int_{\mathbb{S}^{2n+1}} |f|^2 d\xi,
\end{equation*}
where $\mathcal{A}_d(0<d<Q)$ are the intertwining operator, $\mathcal{A}'_Q$ is the conditional intertwinor introduced in \cite{BFM2013}, and $d\xi$ is the normalized surface measure of $\mathbb{S}^{2n+1}$.
\end{abs}\\
\begin{key words} Sharp Sobolev inequality, sharp Hardy-Littlewood-Sobolev inequality, complex sphere, CR manifold
\end{key words}
\section{Introduction}\label{Sec Intro}

It is well known that the classical Sobolev inequalities and Hardy-Littlewood-Sobolev(HLS) inequalities are basic tools in analysis and geometry and their sharp constants play an essential role because they contain geometric and probabilistic information (see e.g., \cite{Be1993, CY1987, T1968, Y1960}).
Recently, many interesting and challenging results on Riemannian geometry and sub-Riemannian manifolds ( such as Heisenberg Group,CR sphere) were  also obtained to understand  different  geometry framework. In particular, many interesting geometric inequalities,
Sobolev-type inequalities and HLS inequality on the sub-Riemannian manifolds attracted the attention of analysts (see e.g., \cite{BFM2013,Folland-Stein1974,Frank-Lieb2012a,CLZ2016,HLZ2012}). Based on  the work of Frank and Lieb \cite{Frank-Lieb2012a} this paper establishes the CR-sphere counterpart of the Sobolev inequalities discussed in \cite{Be1993} in the Euclidean-sphere setting.

For convenience, we firstly introduce some notations and known facts about the complex sphere $\mathbb{S}^{2n+1}$. More details can be found in \cite{BFM2013} and references therein.

Denoted by $\mathbb{S}^{2n+1}$ the complex sphere
    $$\mathbb{S}^{2n+1}= \Bigl\{ \xi=(\xi_1,\xi_2, \cdots,\xi_{n+1}) \in \mathbb{C}^{n+1}: \ \sum_{j=1}^{n+1}|\xi_j|^2=1 \Bigr\}.$$
Then $\mathbb{C}T\mathbb{S}^{2n+1}$ is generated by the vectors $T_j,\overline{T}_j,\ j=1,2,\cdots,n+1$ and $\mathcal{T}$, where
    $$T_j= \frac\partial{\partial\xi_j} -\overline{\xi_j} \sum_{j=1}^{n+1} \xi_k\frac\partial{\partial \xi_k},\ j=1,2,\cdots,n+1, \text{ and } \mathcal{T}= \frac i2 \sum_{k=1}^{n+1} \Bigl(\xi_k\frac\partial{\partial \xi_k}-\overline{\xi_k}\frac\partial{\partial \overline{\xi_k}}\Bigr).$$
Let $Q=2n+2$ be the homogeneous dimension induced from Heisenberg group by Cayley transformation and denote by $d\xi$ the normalized surface measure on $\mathbb{S}^{2n+1}$.

It is known that $L^2(\mathbb{S}^{2n+1})$ can be decomposed into its $U(n+1)$-irreducible components
\begin{equation}\label{docompose}
    L^2(\mathbb{S}^{2n+1})=\bigoplus_{j,k\geq 0}\mathcal{H}_{jk},
\end{equation}
where $\mathcal{H}_{jk}$ is the space of restrictions to $\mathbb{S}^{2n+1}$ of harmonic polynomials $p(z,\bar{z})$ on $\mathbb{C}^{n+1}$ which are homogeneous of degree $j$ in $z$ and degree $k$ in $\bar{z}$. Take $\{Y_{jk}\}$ as an orthonormal basis of $\mathcal{H}_{jk}$. Moreover, denote the {\it Hardy spaces} as follows:
\begin{align*}
    \mathcal{H} &= \bigoplus_{j\geq 0} \mathcal{H}_{j0}\\
    &=\{ L^2\ \text{boundary values of holomorphic functions on the unit ball}\},\\
    \overline{\mathcal{H}} &= \bigoplus_{j\geq 0} \mathcal{H}_{0j}\\
    &=\{L^2\ \text{boundary values of antiholomorphic functions on the unit ball}\},\\
    \mathcal{P} &= \bigoplus_{j>0} (\mathcal{H}_{j0} \oplus \mathcal{H}_{0j}) \bigoplus \mathcal{H}_{00} =\{L^2\ \text{CR-pluriharmonic functions}\},\\
    \mathbb{R}\mathcal{P} &=\{L^2\ \text{real-valued CR pluriharmonic functions}\}.
\end{align*}

For $0<d<Q$, the general {\it intertwining operator} $\mathcal{A}_d$ of order $d$ is defined with respect to the spherical harmonics as
\begin{equation}\label{eigenvalue of intertwining operator}
    \mathcal{A}_d Y_{j,k}=\lambda_j(d)\lambda_k(d) Y_{j,k},\quad j,k=0,1,2,\cdots,
\end{equation}
where
    $$\lambda_j(d)= \frac{\Gamma((Q+d)/4+j)} {\Gamma((Q-d)/4+j)},\quad j=0,1,2,\cdots.$$
In particular, $\mathcal{A}_2$ is the \textsl{conformal sublaplacian} $\mathcal{D} =\mathcal{L}+\frac{n^2}4 =\mathcal{L}+(\lambda_0(2))^2$ with
    $$\mathcal{L}=-\frac 12 \sum_{j=1}^{n+1} (T_j\overline{T}_j+\overline{T}_j T_j).$$
Recently, Branson et al \cite{BFM2013} introduced a class of intertwinors $\mathcal{A}'_Q$ of order $Q$, named {\it conditional intertwinors} and  defined on $\mathcal{P}$ as
\begin{equation}\label{eigenvalue of conditional intertwinor}
    \mathcal{A}'_Q Y_{j0}= \lambda_j(Q) Y_{j0}= j(j+1)\cdots (j+n) Y_{j0},\quad \mathcal{A}'_Q Y_{0k}=\lambda_k(Q) Y_{0k}.
\end{equation}

In \cite{JL1987} and \cite{Frank-Lieb2012a}, two classes of Sobolev inequalities (see Theorem 3.1 and Corollary 2.3 of \cite{Frank-Lieb2012a}) were established as follows:
\begin{gather}
    \mathcal{E}[u]
    \geq \frac{n^2}4 \left( \int_{\mathbb{S}^{2n+1}} |u|^{2Q/(Q-2)} d\xi\right)^{(Q-2)/Q} \label{Jerision-Lee's Sobolev}
\intertext{and}
    \frac{4(q-2)}{Q-2} \mathcal{E}_0[u] +\int_{\mathbb{S}^{2n+1}} |u|^2 d\xi
    \geq \left( \int_{\mathbb{S}^{2n+1}} |u|^{q} d\xi\right)^{2/q},\quad 2< q<\frac{2Q}{Q-2} \label{sub Sobolev p=2}
\end{gather}
where 
$\mathcal{E}[u]=\mathcal{E}_0[u]+\frac{n^2}4 u$ and
    $$\mathcal{E}_0[u]=\frac 12\sum_{j=1}^{n+1} (|T_ju|^2 +|\overline{T}_ju|^2).$$
If we adopt the notations of intertwining operator, inequalities \eqref{Jerision-Lee's Sobolev} and \eqref{sub Sobolev p=2} can be rewrote as:
\begin{gather}
    \int_{\mathbb{S}^{2n+1}} u\mathcal{D}u d\xi \geq \frac{n^2}4 \left( \int_{\mathbb{S}^{2n+1}} |u|^{2Q/(Q-2)} d\xi\right)^{(Q-2)/Q}\label{Sobolev 1}
\intertext{and, for $2< q<\frac{2Q}{Q-2}$,}
     \frac{4(q-2)}{Q-2} \int_{\mathbb{S}^{2n+1}} u\mathcal{L}u d\xi +\int_{\mathbb{S}^{2n+1}} |u|^2 d\xi
    \geq \left( \int_{\mathbb{S}^{2n+1}} |u|^{q} d\xi\right)^{2/q},\label{Sobolev 2}
\end{gather}respectively.

\textbf{What is the Sobolev inequality corresponding to the general intertwining operator $\mathcal{A}_d$?}

To answer this question and motivated by the idea "fractional integration controls Sobolev inequality", we establish firstly the following HLS inequalities.

\begin{theorem}[Subcritical HLS inequalities]\label{thm sub HLS}
Let $0<\lambda<Q=2n+2$ and $\frac{2Q}{2Q-\lambda}<p\leq 2$. Then for any $f,g\in L^p(\mathbb{S}^{2n+1})$, it holds
\begin{equation}\label{sub HLS ineq}
\left|\int_{\mathbb{S}^{2n+1}}\int_{\mathbb{S}^{2n+1}} \frac{\overline{f(\xi)}g(\eta)} {|1-\xi\cdot\bar{\eta}|^{\lambda/2}} d\xi d\eta\right|\leq C_{\lambda,n}\|h\|_p\|g\|_p,
\end{equation}
where
    $$C_{\lambda,n}= \int_{\mathbb{S}^{2n+1}} |1-\xi\cdot\bar{\eta}|^{-\lambda/2} d\eta = \frac{\Gamma(Q/2)\Gamma((Q-\lambda)/2)} {\Gamma^2((2Q-\lambda)/4)}.$$
Moreover, Equality in \eqref{sub HLS ineq} holds if and only if $f$ and $g$ are all constants.
\end{theorem}
\begin{remark}\label{remark HLS}
When $p=\frac{2Q}{2Q-\lambda}$ for $0<\lambda<Q$, then \eqref{sub HLS ineq} is the classical HLS inequalities
\begin{equation}\label{HLS ineq}
    \left|\int_{\mathbb{S}^{2n+1}}\int_{\mathbb{S}^{2n+1}} \frac{\overline{f(\xi)}g(\eta)} {|1-\xi\cdot\overline{\eta}|^{\lambda/2}} d\xi d\eta\right|\leq C_{\lambda,n}\|h\|_p\|g\|_p.
\end{equation}
Moreover, by Theorem 2.2 of \cite{Frank-Lieb2012a}, we know that equality in \eqref{HLS ineq} holds if and only if
\begin{equation}\label{extremal function for comformal case}
    f(\xi)=\frac{c}{|1-\overline{\zeta}\cdot\xi|^{(2Q-\lambda)/2}},\quad g(\eta)=\frac{c'}{|1-\overline{\zeta}\cdot\xi|^{(2Q-\lambda)/2}}
\end{equation}
for some $c,c'\in\mathbb{C}$ and some $\zeta\in\mathbb{C}^{n+1}$ with $|\zeta|<1$ (unless $f\equiv 0$ or $g\equiv 0$).
\end{remark}

Take $f=g=\sum_{j,k\geq 0}Y_{j,k}$ in \eqref{sub HLS ineq} and \eqref{HLS ineq}. Then, we have by \eqref{expansion of HLS} that
\begin{equation}\label{HLS expansion}
    \sum_{j,k\geq 0}\gamma_{j,k}^\lambda \int_{\mathbb{S}^{2n+1}} |Y_{j,k}|^2 d\xi\leq \|f\|_p^2,\quad \frac{2Q}{2Q-\lambda}\leq p\leq 2.
\end{equation}
By a duality argument and letting $\lambda=Q-d$, we get the following Sobolev inequalities on the $\mathbb{S}^{2n+1}$:
\begin{align}\label{Sobolev 3}
    \|f\|_q^2 &\leq \sum_{j,k\geq 0} \frac 1{\gamma_{j,k}^\lambda} \int_{\mathbb{S}^{2n+1}} |Y_{j,k}|^2 d\xi\nonumber\\
    &=\sum_{j,k\geq 0} \frac{\Gamma(j+(2Q-\lambda)/4) \Gamma(k+(2Q-\lambda)/4) \Gamma^2(\lambda/4)} {\Gamma^2((2Q-\lambda)/4) \Gamma(j+\lambda/4) \Gamma(k+\lambda/4)} \int_{\mathbb{S}^{2n+1}} |Y_{j,k}|^2 d\xi\nonumber\\
    &=\sum_{j,k\geq 0} \frac{\Gamma(j+(Q+d)/4) \Gamma(k+(Q+d)/4) \Gamma^2((Q-d)/4)} {\Gamma(j+(Q-d)/4) \Gamma(k+(Q-d)/4) \Gamma^2((Q+d)/4)} \int_{\mathbb{S}^{2n+1}} |Y_{j,k}|^2 d\xi\nonumber\\
    &=\frac 1{(\lambda_0(d))^2} \int_{\mathbb{S}^{2n+1}} f\mathcal{A}_d f d\xi,\quad 2\leq q\leq\frac{2Q}{Q-d}.
\end{align}
Particularly, if $d=2$ and $q=\frac{2Q}{Q-2}$, then \eqref{Sobolev 3} is Sobolev inequality \eqref{Sobolev 1}. While for $d=2$ and $2<q<\frac{2Q}{Q-2}$, we find that the constant $\frac 1{(\lambda_0(2))^2}$ is strictly bigger than the constant $\frac{4(q-2)}{Q-2}$ of \eqref{Sobolev 2} and therefore not sharp. Next theorem gives the sharp form of the Sobolev inequalities on the CR-sphere.

\begin{theorem}\label{thm Sobolev}
For any $f\in C^\infty(\mathbb{S}^{2n+1})$ and $0<d<Q$, we have:

1) {\bf Conformal Sobolev inequalities}: For $q=\frac{2Q}{Q-d}$,
\begin{equation}\label{Sobolev ineq conformal}
   \|f\|_q^2\leq  \frac{\Gamma^2((Q-d)/4)} {\Gamma^2((Q+d)/4)} \int_{\mathbb{S}^{2n+1}} f\mathcal{A}_df d\xi.
\end{equation}
Moreover, equality holds if and only if
\begin{equation}\label{extremal conformal Sobolev}
    f(\xi)=c|1-\bar\zeta\cdot\xi|^{(d-Q)/2}
\end{equation}
for some $c\in\mathbb{C}$ and some $\zeta\in\mathbb{C}^{n+1}$ with $|\zeta|<1$.

2) {\bf Subcritical Sobolev inequalities}: For $2\leq q< \frac{2Q}{Q-d}$,
\begin{align}\label{Sobolev ineq}
    \|f\|_q^2\leq & \frac{8(q-2)}{d(Q-d)} \frac{\Gamma^2((Q-d)/4+1)} {\Gamma^2((Q+d)/4)}\left( \int_{\mathbb{S}^{2n+1}} f\mathcal{A}_df d\xi\right. \nonumber\\
    &\left. -(\lambda_0(d))^2 \int_{\mathbb{S}^{2n+1}} |f|^2 d\xi\right) +\int_{\mathbb{S}^{2n+1}} |f|^2 d\xi.
\end{align}
Moreover, for $2<q<\frac{2Q}{Q-d}$, equality holds if and only if $f$ is constant.
\end{theorem}

\begin{remark}
The conformal Sobolev inequalities \eqref{Sobolev ineq conformal} and their derivation from Frank and Lieb HLS inequality on the Heisenberg group (\cite{Frank-Lieb2012a}) are well known within the group of researchers interested in conformal geometry (see \cite{BFM2013} for further  details). We provide concise proof for completeness. On the other hand the subcritical Sobolev inequalities \eqref{Sobolev ineq} are new. Their Euclidean counterpart can be found in \cite{Be1993}.
\end{remark}

\begin{remark}
When $d=2$, \eqref{Sobolev ineq conformal} and \eqref{Sobolev ineq} are \eqref{Sobolev 1} and \eqref{Sobolev 2}, respectively.
\end{remark}

Combining the method of Beckner in \cite{Be1993} with the HLS inequality on the Heisenberg group (\cite{Frank-Lieb2012a}) and letting $d\rightarrow Q^-$, we have the following sharp inequalities.

\begin{theorem}\label{thm end case}
For any $f\in C^\infty(\mathbb{S}^{2n+1}) \cap\mathbb{R}\mathcal{P}$, we have:

1) {\bf Beckner-Onofri's inequality}:
\begin{equation}\label{Beckner-Onofri}
    \frac 1{2(n+1)!} \int_{\mathbb{S}^{2n+1}} f \mathcal{A}'_Q f d\xi +\int_{\mathbb{S}^{2n+1}} f d\xi-\log \int_{\mathbb{S}^{2n+1}} e^f d\xi\geq 0;
\end{equation}

2) {\bf Subcritical Sobolev inequalities}: for $2\leq q< +\infty$,
\begin{equation}\label{Sobolev ineq sub end}
    \|f\|_q^2\leq \frac{q-2}{(n+1)!} \int_{\mathbb{S}^{2n+1}} f \mathcal{A}'_Q f d\xi +\int_{\mathbb{S}^{2n+1}} |f|^2 d\xi.
\end{equation}
\end{theorem}

\begin{remark}
Note that Beckner-Onofri's inequalities \eqref{Beckner-Onofri} is the main result of \cite{BFM2013}. The authors of \cite{BFM2013} were well aware that \eqref{Beckner-Onofri} could be derived from \eqref{Sobolev ineq conformal} and they also say how, but \cite{BFM2013} was made available as a preprint several years before \cite{Frank-Lieb2012a} was published and at the time \eqref{Sobolev ineq conformal} was only a conjecture. So, for conciseness, we omit the proof.
\end{remark}

\begin{remark}
As in \cite{Be1993}, by making the substitution $f\rightarrow 1+\frac 1q f$ in \eqref{Sobolev ineq sub end} and taking the limit $q\rightarrow +\infty$ for bounded $f$, we can obtain \eqref{Beckner-Onofri} again.
\end{remark}

The plan of the paper is as follows. Section \ref{Sec proof} is devoted to the proof of Theorem \ref{thm sub HLS}, Theorem \ref{thm Sobolev} and the subcritical case of Theorem \ref{thm end case}. Our main tools are the Funck-Heck Theorem on the complex sphere and the duality argument. For completeness, in Appendix \ref{sec decompose}, we state the Fuck-Heck theorem established by Frank and Lieb in \cite{Frank-Lieb2012a} and give some applications.

\section{Proofs of Theorem \ref{thm sub HLS}, Theorem \ref{thm Sobolev} and Theorem \ref{thm end case}}\label{Sec proof}

\noindent {\bf Proof of Theorem \ref{thm sub HLS}.} 1) {\bf Case $\frac{2Q}{2Q-\lambda}<p<2$.}

Firstly, we claim that, for any $\lambda_1$ and $\lambda_2$ satisfying $0<\lambda_1<\lambda_2<Q$ and any $f\in L^2(\mathbb{S}^{2n+1})$, it holds
\begin{equation}\label{formula 3.1}
   \frac{\int_{\mathbb{S}^{2n+1}}\int_{\mathbb{S}^{2n+1}} \frac{\overline{f(\xi)} f(\eta)} {|1-\xi\cdot\bar\eta|^{\lambda_1/2}} d\xi d\eta} {\int_{\mathbb{S}^{2n+1}} |1-\xi\cdot\bar\eta|^{-\lambda_1/2}d\eta} \leq \frac{\int_{\mathbb{S}^{2n+1}}\int_{\mathbb{S}^{2n+1}} \frac{\overline{f(\xi)} f(\eta)} {|1-\xi\cdot\bar\eta|^{\lambda_2/2}} d\xi d\eta}{\int_{\mathbb{S}^{2n+1}} |1-\xi\cdot\bar\eta|^{-\lambda_2/2} d\eta}.
\end{equation}
Moreover, equality holds if and only if $f$ is constant.

Now, Taking $\lambda_1=\lambda$ and $\lambda_2=2Q(1-1/q)$ in \eqref{formula 3.1}, noting the positivity of the left side of \eqref{sub HLS ineq} and combining with the classical HLS inequalities \eqref{HLS ineq}, we can complete the proof of Theorem \ref{thm sub HLS} for the case $\frac{2Q}{2Q-\lambda}<q<2$ since $L^2(\mathbb{S}^{2n+1})$ is dense in $L^q(\mathbb{S}^{2n+1})$. Therefore, it is sufficient to prove \eqref{formula 3.1}.

To prove inequality \eqref{formula 3.1}, we only need to show $\gamma_{j,k}^{\lambda_1}\leq\gamma_{j,k}^{\lambda_2},\ j,k=0,1,2,\cdots$ by \eqref{expansion of HLS}.

Obviously, $\gamma_{0,0}^{\lambda_1}=\gamma_{0,0}^{\lambda_2}$. While for $j+k\geq 1$, it is easy to see that
    $$\gamma_{j,k}^\lambda=\frac{\Gamma^2((2Q-\lambda)/4) \Gamma(j+\lambda/4) \Gamma(k+\lambda/4)} {\Gamma(j+(2Q-\lambda)/4) \Gamma(k+(2Q-\lambda)/4) \Gamma^2(\lambda/4)}$$
is strictly increasing with respect to $\lambda$. Therefore, \eqref{formula 3.1} holds. Moreover, by the decomposition of $L^2$ function, we know that equality in \eqref{formula 3.1} holds if and only if $f$ is a constant.

\smallskip
\noindent 2) {\bf Case $q=2$}

Take the spherical harmonic expansion $f(\xi)=\sum_{j,k\geq 0} Y_{j,k}(\xi)$ with $Y_{j,k}\in\mathcal{H}_{j,k}$. Then inequality \eqref{sub HLS ineq} is equivalent to
    $$\sum_{j,k\geq 0} \gamma_{j,k}^\lambda \int_{\mathbb{S}^{2n+1}} |Y_{j,k}(\xi)|^2 d\xi\leq \sum_{j,k\geq 0} \int_{\mathbb{S}^{2n+1}} |Y_{j,k}(\xi)|^2 d\xi.$$
On the other hand, it is easy to obtain that $\gamma_{0,0}^\lambda=1$ and $\gamma_{j,k}^\lambda<1$ for $j+k\geq 1$. So, we complete the proof. \hfill $\Box$

\medskip

\noindent {\bf Proof of Part 1) of Theorem \ref{thm Sobolev}: Conformal Sobolev inequalities.}

By \eqref{HLS expansion}, we know that, for any $g(\xi)=\sum_{j,k\geq 0}Y_{j,k}(\xi)\in C^\infty(\mathbb{S}^{2n+1})$,
\begin{equation}\label{formula 3.8}
    \sum_{j,k\geq 0}\gamma_{j,k}^\lambda \int_{\mathbb{S}^{2n+1}} |Y_{j,k}|^2 d\xi\leq \|g\|_p^2\quad\text{with}\quad  p=\frac{2Q}{2Q-\lambda}.
\end{equation}
So, for any $f(\xi)=\sum_{j,k\geq 0} Z_{j,k}(\xi)\in C^\infty(\mathbb{S}^{2n+1})$,
\begin{align}\label{formula 3.9}
    &\left|\int_{\mathbb{S}^{2n+1}} \overline{f(\xi)} g(\xi) d\xi\right| =\left|\sum_{j,k\geq 0} \int_{\mathbb{S}^{2n+1}} \overline{Z_{j,k}(\xi)} Y_{j,k}(\xi) d\xi\right|\nonumber\\
    \leq &\sqrt{\sum_{j,k\geq 0} \frac 1{\gamma_{j,k}^\lambda} \int_{\mathbb{S}^{2n+1}} |Z_{j,k}|^2 d\xi} \cdot \sqrt{\sum_{j,k\geq 0} \gamma_{j,k}^{\lambda} \int_{\mathbb{S}^{2n+1}} |Y_{j,k}|^2 d\xi}\nonumber\\
    \leq & \|g\|_p \left(\sum_{j,k\geq 0} \frac 1{\gamma_{j,k}^\lambda} \int_{\mathbb{S}^{2n+1}} |Z_{j,k}|^2 d\xi\right)^{\frac 12}= \|g\|_p \left( \frac 1{(\lambda_0(d))^2} \int_{\mathbb{S}^{2n+1}} f\mathcal{A}_d f d\xi \right)^{\frac 12},
\end{align}
where $d=Q-\lambda\in (0,Q)$. Because of the arbitrariness of $g$ and the density, we get
\begin{equation}\label{formula 3.10}
    \|f\|_q^2\leq \frac 1{(\lambda_0(d))^2} \int_{\mathbb{S}^{2n+1}} f\mathcal{A}_d f d\xi
\end{equation}
for any $f\in L^q(\mathbb{S}^{2n+1})$ and $q=\frac{2Q}{Q-d}$.

A direct computation shows that, if $f$ is defined as in \eqref{extremal conformal Sobolev}, then equality in \eqref{formula 3.10} holds. So, the constant $\frac 1{(\lambda_0(d))^2}$ of \eqref{formula 3.10} is sharp.
In the following we discuss the extremal functions.

Assume nonnegative function $f_0\in L^q(\mathbb{S}^{2n+1})$ be an extremal function of \eqref{formula 3.10}, i.e.,
\begin{equation}\label{formula 3.11}
    \|f_0\|_q^2\leq \frac 1{(\lambda_0(d))^2} \int_{\mathbb{S}^{2n+1}} f_0\mathcal{A}_d f_0 d\xi.
\end{equation}
By \eqref{formula 3.9}, we have
\begin{equation}\label{formula 3.12}
    |<f_0,g>|\leq \|f_0\|_q\|g\|_{q'} \quad\text{with}\quad q'=\frac{2Q}{Q+d}.
\end{equation}
It is know that there exists some function $g_0\in L^{q'}(\mathbb{S}^{2n+1})$ such that equality in \eqref{formula 3.12} holds. Using the property of H\"{o}lder inequality, we know that $f_0 =c g_0 ^{\frac{Q-d} {Q+d}}$, where $c$ is some constant. Substituting $f_0$ and $g_0$ into \eqref{formula 3.9}, we find that $g_0$ is an extremal function of \eqref{HLS ineq}. So, the extremal function $f_0$ must have the form \eqref{extremal conformal Sobolev}. \hfill $\Box$

\medskip
\noindent {\bf Proof of Part 2) of Theorem \ref{thm Sobolev}: Subcritical Sobolev inequalities.}

Note that case $q=2$ is trivial. Therefore, we assume $2<q<\frac{2Q}{Q-d}$ in the sequel. If
\begin{align}\label{formula 3.2}
    &\frac{\Gamma^2((Q-d_1)/4)} {\Gamma^2((Q+d_1)/4)} \int_{\mathbb{S}^{2n+1}} f\mathcal{A}_{d_1}f d\xi\nonumber\\
    \leq & \frac{8(q-2)}{d(Q-d)} \frac{\Gamma^2((Q-d)/4+1)} {\Gamma^2((Q+d)/4)}\left( \int_{\mathbb{S}^{2n+1}} f\mathcal{A}_df d\xi\right. \nonumber\\
    &\left. -(\lambda_0(d))^2 \int_{\mathbb{S}^{2n+1}} |f|^2 d\xi\right) +\int_{\mathbb{S}^{2n+1}} |f|^2 d\xi
\end{align}
holds for $d_1=Q(1-2/q)$ and $\frac{2Q}{Q-d}>q>2$, then we can get \eqref{Sobolev ineq} by combining \eqref{Sobolev ineq conformal}. For showing inequality \eqref{formula 3.2}, by the definition of operator $\mathcal{A}_d$, we need to  prove
\begin{equation*}
    \frac{\lambda_j(d_1)\lambda_k(d_1)} {(\lambda_0(d_1))^2} \leq 1+\frac{8(q-2)}{d(Q-d)} \frac{\Gamma^2((Q-d)/4+1)} {\Gamma^2((Q+d)/4)} (\lambda_j(d)\lambda_k(d) -(\lambda_0(d))^2),
\end{equation*}for $j,k\geq 0$,
So, we will prove that, for $j,k\geq 0$,
\begin{align}\label{formula 3.3}
    &\frac{\Gamma(j+\frac Q{2q'}) \Gamma(k+\frac Q{2q'}) \Gamma^2(\frac Q{2q})} {\Gamma(j+\frac Q{2q}) \Gamma(k+\frac Q{2q}) \Gamma^2(\frac Q{2q'})}\nonumber\\
    \leq & 1+\frac{8(q-2)}{d(Q-d)} \frac{\Gamma^2(\frac{Q-d}4+1)} {\Gamma^2(\frac{Q+d}4)} \left(\frac{\Gamma(j+\frac{Q+d}4) \Gamma(k+\frac{Q+d}4)} {\Gamma(j+\frac{Q-d}4) \Gamma(k+\frac{Q-d}4)}-\frac{\Gamma^2(\frac{Q+d}4)} {\Gamma^2(\frac{Q-d}4)} \right),
\end{align}
where $q'$ is the conjugate number of $q$, i.e., $\frac 1q+\frac 1{q'}=1$. A direct calculation shows that equality in \eqref{formula 3.3} occurs at $(j,k)=(0,0),(1,0)$ or $(0,1)$.

To prove \eqref{formula 3.3}, we differentiate with respect to $j$ and $k$. If the left derivation is less than the right for $j+k\geq 1$, then we can deduce \eqref{formula 3.3} for all $j,k\geq 0$ from the monotonicity. In fact,
\begin{align}\label{formula 3.4}
    &\frac{\partial}{\partial k} \left( \frac{\Gamma(j+\frac Q{2q'}) \Gamma(k+\frac Q{2q'}) \Gamma^2(\frac Q{2q})} {\Gamma(j+\frac Q{2q}) \Gamma(k+\frac Q{2q}) \Gamma^2(\frac Q{2q'})} \right)\nonumber\\
    =& \frac{\Gamma(j+\frac Q{2q'}) \Gamma(k+\frac Q{2q'}) \Gamma^2(\frac Q{2q})} {\Gamma(j+\frac Q{2q}) \Gamma(k+\frac Q{2q}) \Gamma^2(\frac Q{2q'})} \left( \frac{\Gamma'(k+\frac Q{2q'})} {\Gamma(k+\frac Q{2q'})} -\frac{\Gamma'(k+\frac Q{2q})} {\Gamma(k+\frac Q{2q})}\right)\nonumber\\
    =& \frac{\Gamma(j+\frac Q{2q'}) \Gamma(k+\frac Q{2q'}) \Gamma^2(\frac Q{2q})} {\Gamma(j+\frac Q{2q}) \Gamma(k+\frac Q{2q}) \Gamma^2(\frac Q{2q'})} \sum_{l=0}^{+\infty}\left( \frac{1}{k+\frac Q{2q}+l}-\frac 1{k+\frac Q{2q'}+l}\right)\nonumber\\
    =& \frac{\Gamma(j+\frac Q{2q'}) \Gamma(k+\frac Q{2q'}) \Gamma^2(\frac Q{2q}) \frac Q{2q}} {\Gamma(j+\frac Q{2q}) \Gamma(k+\frac Q{2q}) \Gamma^2(\frac Q{2q'})}\sum_{l=0}^{+\infty} \frac {q-2} {(l+k)^2+\frac Q2(l+k)+(\frac Q2)^2\frac 1q\frac 1{q'}},
\end{align}
and
\begin{align}\label{formula 3.5}
    &\frac\partial{\partial k}\left[1+\frac{8(q-2)}{d(Q-d)} \frac{\Gamma^2(\frac{Q-d}4+1)} {\Gamma^2(\frac{Q+d}4)} \left(\frac{\Gamma(j+\frac{Q+d}4) \Gamma(k+\frac{Q+d}4)} {\Gamma(j+\frac{Q-d}4) \Gamma(k+\frac{Q-d}4)}-\frac{\Gamma^2(\frac{Q+d}4)} {\Gamma^2(\frac{Q-d}4)} \right)\right]\nonumber\\
    =&\frac{8(q-2)}{d(Q-d)} \frac{\Gamma^2(\frac{Q-d}4+1)} {\Gamma^2(\frac{Q+d}4)} \frac{\Gamma(j+\frac{Q+d}4) \Gamma(k+\frac{Q+d}4)} {\Gamma(j+\frac{Q-d}4) \Gamma(k+\frac{Q-d}4)}\left( \frac{\Gamma'(k+\frac{Q+d}4)} {\Gamma(k+\frac{Q+d}4)} -\frac{\Gamma'(k+\frac{Q-d}4)} {\Gamma(k+\frac{Q-d}4)}\right)\nonumber\\
    =& \frac{\frac{Q-d}4 \Gamma^2(\frac{Q-d}4)} {\Gamma^2(\frac{Q+d}4)} \frac{\Gamma(j+\frac{Q+d}4) \Gamma(k+\frac{Q+d}4)} {\Gamma(j+\frac{Q-d}4) \Gamma(k+\frac{Q-d}4)} \sum_{l=0}^{+\infty} \frac {q-2} {(l+k)^2+\frac Q2(l+k)+ \frac{Q-d}{4} \frac{Q+d}{4}}.
\end{align}
Combining the facts:  $\frac{d}{dx}\frac{\Gamma(l+x)}{\Gamma(x)}\geq 0$ for $x>0$ and $l\geq 0$, $\frac{d}{dx}\frac{\Gamma(l+x)} {\Gamma(1+x)}\geq 0$ for $x>0$ and $l\geq 1$, and $\frac{2Q}{Q+d}<q'<2<q<\frac{2Q}{Q-d}$, we have, for $j+k\geq 1$
\begin{equation}\label{formula 3.6}\begin{cases}
    \frac{\Gamma(j+\frac Q{2q'}) \Gamma(k+\frac Q{2q'})} {\Gamma^2(\frac Q{2q'})}\leq \frac{\Gamma(j+\frac{Q+d}4) \Gamma(k+\frac{Q+d}4)} {\Gamma^2(\frac{Q+d}4)},\\
    \frac{\Gamma(j+\frac Q{2q}) \Gamma(k+\frac Q{2q})} {\frac Q{2q} \Gamma^2(\frac Q{2q})} \geq \frac{\Gamma(j+\frac{Q-d}4) \Gamma(k+\frac{Q-d}4)} {\frac{Q-d}4 \Gamma^2(\frac{Q-d}4)}.
\end{cases}\end{equation}
Moreover, since $f(x)=x(1-x)$ is strictly increasing on $[0,\frac 12]$, then
    $$\frac{Q-d}{2Q} \cdot\frac{Q+d}{2Q} =f(\frac{Q-d}{2Q}) <f(\frac 1q) =\frac 1q \cdot\frac 1{q'},$$
which implies that
\begin{equation}\label{formula 3.7}
    \frac {q-2} {(l+k)^2+\frac Q2(l+k)+(\frac Q2)^2\frac 1q\frac 1{q'}} \leq \frac {q-2} {(l+k)^2+\frac Q2(l+k)+ \frac{Q-d}{4} \frac{Q+d}{4}}
\end{equation}
for $k\geq 0$ and $l\geq 0$. So the $k$ derivative of the LHS of \eqref{formula 3.3} is less of the one of the RHS and, the same is true for the $j$ derivative. Then, we get \eqref{formula 3.3}.

From the above proof, we know that equality of \eqref{formula 3.3} occurs only at $(j,k)=(0,0),\ (1,0)$ or $(0,1)$. Therefore, equality of \eqref{formula 3.2} holds if and only if
    $$f\in\mathcal{H}_{0,0} \bigoplus\mathcal{H}_{0,1} \bigoplus\mathcal{H}_{1,0}.$$
Combining the extremal result of \eqref{extremal conformal Sobolev}, we know that equality of \eqref{Sobolev ineq} for $2<q<\frac{2Q}{Q-d}$ holds if and only if $f$ is constant. \hfill $\Box$

\medskip
\noindent {\bf Proof of part 2) of Theorem \ref{thm end case}: Subcritical Sobolev inequalities.}

For any $Y_{j0}\in\mathcal{H}_{j0}$, $j=0,1,2,\cdots$, we have, as $d\rightarrow Q^-$,
\begin{align*}
    &\frac{8(q-2)}{d(Q-d)} \frac{\Gamma^2((Q-d)/4+1)} {\Gamma^2((Q+d)/4} \mathcal{A}_d Y_{j0}\\ =&\frac{2(q-2)}d \frac{\frac{Q+d}4 (\frac{Q+d}4+1)\cdots (\frac{Q+d}4+j-1)} {(\frac{Q-d}4+1)\cdots (\frac{Q-d}4+j-1)} Y_{j0}\\
    \rightarrow &(q-2)\frac{j(j+1)\cdots(j+n)} {(n+1)!} Y_{j0}= \frac{q-2} {(n+1)!} \mathcal{A}'_Q Y_{j0}.
\end{align*}
Similarly, the above result holds for any $Y_{0k}\in \mathcal{H}_{0k}$, $k=0,1,2,\cdots$. On the other hand, we have
    $$\lambda_0(d)\rightarrow 0, \quad\text{as}\quad d\rightarrow Q^-.$$
So, we get \eqref{Sobolev ineq sub end} via letting $d\rightarrow Q^{-}$ in \eqref{Sobolev ineq}.

\begin{appendix}
\section{The Funk-Hecke Theorem on the complex sphere}\label{sec decompose}

In \cite{Frank-Lieb2012a}, Frank and Lieb established the following two results. Notice that, in the following formulas, the factor $|\mathbb{S}^{2n+1}|$ appears in the denominators because we use the normalized surface measure.

\begin{proposition}[Proposition 5.2 of \cite{Frank-Lieb2012a}]
Let $K$ be an integrable function on the unit ball in $\mathbb{C}$. Then the operator on $\mathbb{S}^{2n+1}$ with kernel $K(\xi\cdot\overline{\eta})$ is diagonal with respect to decomposition \eqref{docompose}, and on the space $\mathcal{H}_{j,k}$ its eigenvalue is given by
\begin{equation}\label{eigenvalue 1}\begin{split}
    \frac 1{|\mathbb{S}^{2n+1}|}\frac{\pi^n m!}{2^{n+|j-k|/2}(m+n-1)!} &\int_{-1}^1 dt (1-t)^{n-1} (1+t)^{|j-k|/2}P_m^{(n-1,|j-k|)}(t)\\
    &\times\int_{-\pi}^\pi d\varphi K(e^{-i\varphi}\sqrt{(1+t)/2})e^{i(j-k)\varphi},
\end{split}\end{equation}
where $m:=\min\{j,k\}$ and $P_m^{(\alpha,\beta)}$ are the Jacobi polynomials.
\end{proposition}

\begin{proposition}[Corollary 5.3 of \cite{Frank-Lieb2012a}]\label{thm specific eigenvalue}
Let $-1<\alpha<\frac{n+1}2$.

(1) The eigenvalue of the operator with kernel $|1-\xi\cdot\overline{\eta}|^{-2\alpha}$ on the subspace $\mathcal{H}_{j,k}$ is
\begin{equation}\label{eigenvalue 3}
    E_{j,k}:=\frac{2\pi^{n+1}\Gamma(n+1-2\alpha)} {|\mathbb{S}^{2n+1}|\Gamma^2(\alpha)} \frac{\Gamma(j+\alpha)}{\Gamma(j+n+1-\alpha)} \frac{\Gamma(k+\alpha)}{\Gamma(k+n+1-\alpha)}.
\end{equation}

(2)The eigenvalue of the operator with kernel $|\xi\cdot\overline{\eta}|^2|1-\xi\cdot\overline{\eta}|^{-2\alpha}$ on the subspace $\mathcal{H}_{j,k}$ is
\begin{equation}\label{eigenvalue 2}
    E_{j,k}\left(1- \frac{(\alpha-1)(n+1-2\alpha)(2jk+n(j+k-1+\alpha)} {(j-1+\alpha)(j+n+1-\alpha)(k-1+\alpha)(k+n+1-\alpha)}\right).
\end{equation}
When $\alpha=0$ or $1$, formula \eqref{eigenvalue 2} and \eqref{eigenvalue 3} are to be understood by taking limits with fixed $j$ and $k$.
\end{proposition}

As application, we have the following result.

\begin{proposition}\label{thm constants}
For $0<\lambda<Q$, we have
\begin{equation}\label{constant HLS}
    \int_{\mathbb{S}^{2n+1}} |1-\xi\cdot\bar\eta|^{-\lambda/2} d\eta = \frac{\Gamma(Q/2)\Gamma((Q-\lambda)/2)} {\Gamma^2((2Q-\lambda)/4)}.
\end{equation}
For $f(\xi)=\sum_{j,k\geq 0} Y_{j,k}$ with $Y_{j,k}\in\mathcal{H}_{j,k}$, then
\begin{equation}\label{expansion of HLS}
    \frac{\int_{\mathbb{S}^{2n+1}}\int_{\mathbb{S}^{2n+1}} \frac{{f(\xi)}f(\eta)} {|1-\xi\cdot\bar\eta|^{\lambda/2}} d\xi d\eta}{\int_{\mathbb{S}^{2n+1}} |1-\xi\cdot\bar\eta|^{-\lambda/2} d\eta} =\sum_{j,k\geq 0} \gamma_{j,k}^\lambda \int_{\mathbb{S}^{2n+1}} |Y_{j,k}(\xi)|^2 d\xi
\end{equation}
with
\begin{equation*}
    \gamma_{j,k}^\lambda=\frac{\Gamma^2((2Q-\lambda)/4) \Gamma(j+\lambda/4) \Gamma(k+\lambda/4)} {\Gamma(j+(2Q-\lambda)/4) \Gamma(k+(2Q-\lambda)/4) \Gamma^2(\lambda/4)},\quad j,k=0,1,2,\cdots.
\end{equation*}
\end{proposition}
\end{appendix}

\medskip

\noindent {\bf Acknowledgements}\\
\noindent
The project is supported by  the
National Natural Science Foundation of China (Grant No. 11201443) and Natural Science Foundation of Zhejiang Province(Grant No. LY18A010013). We would like to thank Professor Meijun Zhu and Professor Jingbo Dou for some helpful discussions. We also thank the referee for his/her careful reading of the original manuscript.


\begin{thebibliography}{a}


\bibitem{Be1993}W. Beckner, Sharp Sobolev inequalities on the sphere and the Moser-Trudinger inequality, Annals of Mathematics, 138(1993), 213-242.

\bibitem{BFM2013}T.P. Branson, L. Fontana and C. Morpurgo, Moser-Trudinger and Beckner-Onofri's inequalities on the CR sphere, Annals of Mathematics, 177(2013), 1-52.

\bibitem{CY1987}S-Y A. Chang and P. Yang, Prescribing Gaussian curvature on $\mathbb{S}^2$, Acta Mathematica, 159(1987), 215-259.


\bibitem{Folland 1975}G.B. Folland, Spherical harmonic expansion of the Possion-Szeg\H{o} kernel for the ball, Proc. Amer. Math. Soc., 47(2)(1975), 401-408.

\bibitem{Folland-Stein1974}G.B. Folland and E.M. Stein, \textsl{Estimates for the $\bar{\partial}_b$ complex and analysis on the Heisenberg group}, Communications on Pure and applied Mathematics, 27(1994), 429-522.

\bibitem{Frank-Lieb2012a}R.L. Frank and E.H. Lieb, \textsl{Sharp constants in several inequalities on the Heisenberg group}, Annals of Mathematics, 176(2012), 349-381.

\bibitem{Frank-Lieb2012b}R.L. Frank and E.H. Lieb, A new, rearrangement-free proof of the sharp Hardy-Littlewood-Sobolev inequality, Operator Theory: Advances and Applications, 219(2012), 55-67.

\bibitem{CLZ2016} M. Christ, H.  Liu,  A. Zhang, Sharp Hardy-Littlewood-Sobolev inequalities on the
octonionic Heisenberg group, Calc. Var. PDE, 55(11)(2016),1-18.
\bibitem{HLZ2012} X. Han, G. Lu, J. Zhu, Hardy-Littlewood-Sobolev and Stein-Weiss inequalities and integral systems on
the Heisenberg group, Nonlinear Anal. 75 (11) (2012), 4296-4314.



\bibitem{JL1987} D. Jerison, J. M. Lee, The Yamabe problem on CR manifolds, J. Differential Geom. 25 (1987): 167-197.

\bibitem{JL1988} D. Jerison, J. M. Lee, Extremals for the Sobolev inequality on the Heisenberg group and the CR Yamabe problem, J. Amer. Math. Soc. 1 (1988): 1-13.

\bibitem{Lee-Parker}J.M. Lee and T.H. Parker, The Yamabe problem, Bull. Amer. Math. Soc. (N.S.), 17(1987): 37-91.

\bibitem{Lieb1983} E. H. Lieb, Sharp constants in the Hardy-Littlewood-Sobolev and related inequalities, Ann. of Math. 118 (1983), 349-374.




\bibitem{T1968} N. S. Trudinger, Remarks concerning the conformal deformation of Riemannian structures on compact manifolds. Annali Sc. Norm. supp. Pisa 22 (1968), 265-274.

\bibitem{Y1960} H. Yamabe, On the deformation of Riemannian structures on compact manifolds, Osaka Math. J., 12 (1960), 21-37.

\end{thebibliography}
\end{document}